\documentclass[11pt]{article}
\usepackage{bez123,calc,curves,ebezier,epic,eepic,graphicx,multiply,rotating}
\usepackage{algorithmic}
\usepackage{bez123,calc,curves,ebezier,epic,eepic,graphicx,multiply,rotating}
\everymath{\displaystyle}
\usepackage{graphicx}
\usepackage{pst-all}
\usepackage{color}
\usepackage{amssymb}
\usepackage{amsmath}
\newtheorem{prethm}{{\bf Theorem}}

\newenvironment{thm}{\begin{prethm}{\hspace{-0.5
               em}{\bf .}}}{\end{prethm}}
\newtheorem{prelemma}{{\bf Lemma}}

\newenvironment{lemma}{\begin{prelemma}{\hspace{-0.5
               em}{\bf .}}}{\end{prelemma}}
\newtheorem{preex}{{\bf Example}}

\newtheorem{preprop}{{\bf Proposition}}

\newtheorem{precor}{{\bf Corollary}}

\newtheorem{preremark}{{\bf Remark}}

\newenvironment{remark}{\begin{preremark}{\hspace{-0.5
               em}{\bf.}}}{\end{preremark}}
\newtheorem{preprob}{{\bf Problem}}

\newtheorem{predefin}{{\bf Definition}}

\newenvironment{defin}{\begin{predefin}{\hspace{-0.5
               em}{\bf .}}}{\end{predefin}}
\newtheorem{preconj}{{\bf Conjecture}}

\newtheorem{preprobb}{{\bf Problem}}

\newtheorem{prelem}{{\bf Theorem}}

\newenvironment{proof}{{\bf Proof.}\rm }{\hfill{$\Box$}}

\newtheorem{presolution}{{\bf Solution.}}

\def\newpic#1{}
\def\qed{\ifhmode\unskip\nobreak\fi\quad\ifmmode\Box\else$\Box$\fi}

\title{\Large\bf\noindent Partial vertex covers and the complexity of some problems concerning static and dynamic monopolies}

\author{\large\bf Hossein Soltani~~~~ Manouchehr Zaker\footnote{mzaker@iasbs.ac.ir}
\vspace{5mm}\\
    Department of Mathematics,\\
     Institute for Advanced Studies in Basic Sciences,\\
    Zanjan 45137-66731, Iran\\
  }
    \date{}
\begin{document}
\maketitle

\begin{abstract}
\noindent Let $G$ be a graph and $\tau$ be an assignment of nonnegative integer thresholds to the vertices of $G$. Denote the average of thresholds in $\tau$ by $\bar{\tau}$. A subset of vertices $D$ is said to be a $\tau$-dynamic monopoly, if $V(G)$ can be partitioned into subsets $D_0, D_1, \ldots, D_k$ such that $D_0=D$ and for any $i\in \{0, \ldots, k-1\}$, each vertex $v$ in $D_{i+1}$ has at least $\tau(v)$ neighbors in $D_0\cup \ldots \cup D_i$. Denote the size of smallest $\tau$-dynamic monopoly by $dyn_{\tau}(G)$. Also a subset of vertices $M$ is said to be a $\tau$-static monopoly (or simply $\tau$-monopoly) if any vertex $v\in V(G)\setminus M$ has at least $\tau(v)$ neighbors in $M$. Denote the size of smallest $\tau$-monopoly by $mon_{\tau}(G)$. For a given positive number $t$, denote by $Sdyn_t(G)$ (resp. $Smon_t(G)$), the minimum $dyn_{\tau}(G)$ (resp. $mon_{\tau}(G)$) among all threshold assignments $\tau$ with $\overline{\tau}\geq t$. In this paper we consider the concept of partial vertex cover as follows. Let $G=(V, E)$ be a graph and $t$ be any positive integer. A subset $S\subseteq V$ is said to be a $t$-partial vertex cover of $G$, if $S$ covers at least $t$ edges of $G$. Denote the smallest size of a $t$-partial vertex cover of $G$ by $P\beta_t(G)$.
Let $\rho$, $0<\rho<1$ be any fixed number and $G$ be a given bipartite graph with $m$ edges. We first prove that to determine the smallest cardinality of a set $S\subseteq V(G)$ such that $S$ covers at least $\rho m$ edges of $G$, is an NP-hard problem. Then we prove that for any constant $t$, $Sdyn_{t}(G)=P\beta_{nt-m}(G)$ and $Smon_t(G)=P\beta_{nt/2}(G)$, where $n$ and $m$ are the order and size of $G$, respectively. Let $\epsilon(G)$ be the edge density of $G$. We finally show that to determine $Sdyn_{k\epsilon(G)}(G)$ ($1<k<2$) and $Smon_{k\epsilon(G)}(G)$ ($0<k<2$) is NP-hard for various classes of graphs, but polynomial-time for trees.
\end{abstract}




\section{Introduction}

\noindent All graphs in this paper are simple and undirected graph. For any graph $G=(V,E)$ we denote the order of $G$ by $|G|$ and the edge density of $G$ by $\epsilon(G):=|E|/|G|$. For other graph theoretical notations, not defined in this paper we refer the reader to \cite{BM}. Throughout the paper by $(G,\tau)$ we mean any graph $G$ together with an assignment of thresholds $\tau: V(G) \rightarrow \Bbb{N}\cup \{0\}$ to the vertices of $G$ such that $0\leq\tau(v)\leq d_{G}(v)$, where $d_{G}(v)$ denotes the degree of vertex $v$ in $G$. Also by the average threshold of $\tau$, denoted by $\overline{\tau}$, we mean ${\sum}_{v\in V(G)} \tau(v)/|G|$. A subset of vertices $D$ is said to be a {\it $\tau$-dynamic monopoly} of $G$, if for some nonnegative integer $k$, the vertices of $G$ can be partitioned into subsets $D_0, D_1, \ldots, D_k$ such that $D_0=D$ and for any $i$, $1\leq i \leq k$, the set $D_{i}$ consists of all vertices $v$ which has at least $\tau(v)$ neighbors in $D_0\cup \ldots \cup D_{i-1}$. Denote the size of minimum $\tau$-dynamic monopoly of $G$ by $dyn_{\tau}(G)$. Dynamic monopolies are in fact modeling the spread of influence in social networks. The spread of innovation or a new product in a community, spread of opinion in Yes-No elections, spread of virus in the internet, spread of disease among population and etc are some examples of these phenomena. Dynamic monopolies and the equivalent concepts target set selection and irreversible iterated process were widely studied in the literature specially in recent years (see e.g. \cite{ABW,CL,CDPRS,CR,C,DR,FKRRS,KSZ,NNUW,Z,Z2}). Also the static version of dynamic monopolies were already known and studied in the literature under the name of monopolies or stamo, of course for simple or strict majority thresholds (see e.g. \cite{FKRRS,KNSZ}). In the following we consider monopolies with general threshold assignments. Given $(G,\tau)$, by a $\tau$-static monopoly (or simply $\tau$-monopoly) in $(G,\tau)$ we mean any subset of vertices $M$ such that for any $v\in V(G)\setminus M$, $d_M(v)\geq \tau(v)$.
Dynamic monopolies with given average threshold were studied in \cite{KSZ, Z2, CR}. Specially in \cite{Z2}, the notation $dyn_t(G)=\min_{\tau: \bar{\tau}=t} dyn_{\tau}(G)$ was introduced, where $t$ is any fixed non-negative number such that $t|G|$ is integer. The following result was proved in \cite{Z2}, where $n$ and $m$ are the order and size of $G$, respectively and $H  \unlhd G$ means $H$ is an induced subgraph of $G$.
$$\min_{\tau: \bar{\tau}=t} dyn_{\tau}(G)= n- \max \{|H|:H  \unlhd G,~|E(H)|\leq 2m-nt\}.$$

\noindent In the definition of $dyn_t(G)$, $t|G|$ is integer. In order to consider all values of $t$, we modify a little bit the definition. But we are forced to make a new notation, i.e. $Sdyn_t(G)$ (which stands for the smallest dynamic monopoly). The formal definition is as follows.

\begin{defin}
Let $G$ be a graph and $t$ a positive number. We define $Sdyn_t(G) = \min\{dyn_{\tau}(G)| \overline{\tau}\geq t\}$. Assume that a subset $D\subseteq V(G)$ and an assignment of thresholds $\tau_0$ with $\overline{\tau_0}\geq t$ are such that $|D|=dyn_{\tau_0}(G)=Sdyn_t(G)$ and $D$ is a $\tau_0$-dynamic monopoly of $(G,\tau_0)$. Then we say $(D, \tau_0)$ is a $t$-Sdynamo of $G$.
\end{defin}

\noindent We have the analogous notation for monopolies.

\begin{defin}
Let $G$ be a graph and $t$ a positive number. We define $Smon_t(G) = \min\{mon_{\tau}(G)| \overline{\tau}\geq t\}$. Assume that a subset $M\subseteq V(G)$ and an assignment of thresholds $\tau_0$ with $\overline{\tau_0}\geq t$ are such that $|M|=mon_{\tau_0}(G)=Smon_t(G)$ and $M$ is a $\tau_0$-monopoly of $(G,\tau_0)$. Then we say $(M, \tau_0)$ is a $t$-Smon of $G$.
\end{defin}

\noindent {\bf The paper is organized as follows.} Section 2 is devoted to partial vertex covers in graphs. We show that to determine the smallest cardinality of a subset which covers a constant fraction of the number of edges in bipartite graphs is NP-hard. The same problem is NP-hard even for the class of planar and the class of chordal graphs. In Section 3, we obtain a relationship between $Smon_t(G)$ and partial vertex covers in $G$ and then prove some hardness results concerning $Smon_t(G)$. In Section 4, we first derive a relationship for $Sdyn_{t}(G)$ in terms of partial vertex covers and then prove NP-completeness of determining $Sdyn_{k\epsilon}(G)$ for bipartite, planar and chordal graphs.

\section{On partial vertex covers}

\noindent Let $G$ be a graph, we say a vertex $v$ covers an edge $e$ of $G$ if $e$ is incident with $v$. A subset of vertices $S$ in a graph $G$, is said to be a vertex cover of $G$ if any edge of $G$ has at least one endpoint in $S$. To determine the smallest size of a vertex cover is classic NP-complete problem even for planar graphs \cite{GJ}. Let $G$ be a graph and $t$ be a positive integer. By a $t$-partial vertex cover we mean any subset $S\subset V(G)$ such that $S$ covers at least $t$ edges of $G$. Denote the smallest size of a $t$-partial vertex cover of $G$ by $P\beta_t(G)$. The computational complexity of problems concerning partial vertex covers were widely studied in the literature (e.g. \cite{GNW,AS,CS}). In the following by PVCB, we mean a decision problem which, given a bipartite graph $G$ and two integers $k$ and $t$, determines whether $G$ contains a $t$-partial vertex cover with at most $k$ vertices. We denote a typical instance of PVCB by $<G,k,t>$. The NP-completeness of PVCB was proved in \cite{AS,CS}. In this paper we consider another variant of partial vertex cover in which we seek for the minimum size of a set $S$ which covers a fixed fraction of $|E(G)|$ edges. Denote this problem by PVC($\rho$), where $\rho$ is any fixed number with $0<\rho<1$.

\noindent {\bf Name:} $PVC(\rho)$\\
\noindent {\bf Instance:} A graph $H=(V,E)$ and a positive integer $l$.\\
\noindent {\bf Question:} Is there a subset $C\subseteq V$ , $|C|\leq l$ such that $C$ covers at least $\rho |E|$ edges of $H$?

\noindent The aim of this section is to prove Theorem \ref{PVCRHO} and Remark \ref{rem1}. In the next sections we obtain some relationships between partial vertex covers and monopolies and also dynamic monopolies of graphs with a given average threshold. We begin with the following lemma.

\begin{lemma}\label{shakh}
Let $G$ be a graph and $G'$ be obtained from $G$ as follows. Corresponding to each vertex $v$ of $G$, we add three vertices $v^1, v^2, v^3$ and edges $vv^1, vv^2, vv^3$ to the graph $G$. Then $<G', k, t+3k>$ is a yes-instance of PVCB if and only if $<G, k, t>$ is a yes-instance of PVCB.
\end{lemma}

\noindent \begin{proof}
First let $<G, k, t>$ be a yes-instance of PVCB. Then there exists a subset $C\subseteq V(G)$ with $|C|=k$ such that $C$ covers at least $t$ edges of $G$. By the construction of $G'$, obviously $C$ covers at least $t+3k$ edges of $G'$. This shows that $<G', k, t+3k>$ is a yes-instance of PVCB.
Now let $<G', k, t+3k>$ be a yes-instance of PVCB. Let also $C \subseteq V(G')$ with  $|C|=k$ covers at least $t+3k$ edges of $G'$. Without loss of generality we may assume that $C \subseteq V(G)$, since we can replace any vertex in $C$ of type $v^j$, $1\leq j\leq 3$, by $v$. It is easily seen that $C$ covers at least $t$ edges of $G$. Hence $<G, k, t>$ is a yes-instance for PVCB.
\end{proof}

\noindent We also need the following lemma, where $\rho$ is any fixed number with $0<\rho<1$ and $\cup$ stands for the vertex disjoint union of graphs.

\begin{lemma}\label{lem2}
Let $G$ be a graph on $n$ vertices with $m$ edges and $G'$ be the graph constructed from $G$ as described in Lemma \ref{shakh}. Let also $t$ and $k$ be two non-negative integers with $k\leq n$. Set $r=\lceil (\rho/(1-\rho))((n(n-1)/2)+3n) \rceil + n+3$ and $s=\lfloor (t+3k+(1-\rho)r+1-\rho(m+3n))/\rho \rfloor$. Let $K_{1,r}$ be a star graph and $v$ be its central vertex. Let also $P_s$ be a path on $s$ vertices and $u$ be a vertex of degree one in $P_s$. Let $w$ be any vertex of degree one in $G'$. The graph $H$ is obtained by adding edges $vu$ and $vw$ to the graph $G'\cup K_{1,r} \cup P_s$ (see Fig. \ref{fig1}). Then $<G', k, t+3k>$ is yes-instance of PVCB if and only if $<H, k+1>$ is yes-instance of PVC($\rho$).
\end{lemma}

\begin{figure}[ht]
\begin{center}
\scalebox{1} 
{
\begin{pspicture}(0,-1.02)(4.604111,1.02)
\psline[linewidth=0.03cm](2.66,0.059999984)(3.38,0.059999984)
\psline[linewidth=0.03cm](4.54,0.059999984)(3.82,0.059999984)
\pscircle[linewidth=0.03,dimen=outer](1.02,0.0){1.02}
\psline[linewidth=0.03cm](2.02,0.06)(2.64,0.06)
\psline[linewidth=0.03cm](2.62,0.06)(3.16,-0.475957)
\psline[linewidth=0.03cm](2.62,0.04)(2.94,-0.575957)
\psline[linewidth=0.03cm](2.62,0.06)(2.3,-0.555957)
\psline[linewidth=0.03cm](2.6,0.02)(2.1,-0.435957)
\psdots[dotsize=0.1](2.94,-0.58000004)
\psdots[dotsize=0.1](3.16,-0.48)
\psdots[dotsize=0.1](2.3,-0.56)
\psdots[dotsize=0.1](2.08,-0.46)
\psdots[dotsize=0.1](2.62,0.04)
\psdots[dotsize=0.1](2.04,0.06)
\psdots[dotsize=0.1](4.54,0.059999984)
\psdots[dotsize=0.1](3.04,0.059999984)
\psdots[dotsize=0.1](4.14,0.059999984)
\usefont{T1}{ptm}{m}{n}
\rput(1.0029101,-0.015){$G'$}
\psdots[dotsize=0.04](2.46,-0.64)
\psdots[dotsize=0.04](2.74,-0.64)
\psdots[dotsize=0.04](3.52,0.059999984)
\psdots[dotsize=0.04](3.62,0.059999984)
\psdots[dotsize=0.04](3.72,0.059999984)
\usefont{T1}{ptm}{m}{n}
\rput(2.5849707,0.26500002){\small $v$}
\usefont{T1}{ptm}{m}{n}
\rput(2.1549706,0.285){\small $w$}
\psdots[dotsize=0.042](2.6,-0.655957)
\usefont{T1}{ptm}{m}{n}
\rput(3.0474854,0.269043){\small $u$}
\end{pspicture}
}
\caption{graph $H$}\label{fig1}

\end{center}
\end{figure}

\noindent \begin{proof}
First note that $s\geq 1$. Let $<G', k, t+3k>$ be a yes-instance of PVCB. Then there exists a subset $C\subseteq V(G')$ with $|C|=k$ such that $C$ covers at least $t+3k$ edges of $G'$. Without loss of generality we may assume that $C$ does not contain $w$ (otherwise, we replace $w$ by the only neighbor of $w$ in $G'$). Obviously $C\cup \{v\}$ covers at least $t+3k+r+2$ edges of $H$. Hence if we show that $t+3k+r+2\geq \lceil \rho|E(H)|\rceil$, then it follows that $<H, k+1>$ is a yes-instance of the PVC($\rho$). For this purpose we show that $t+3k+r+2\geq \rho|E(H)|$, since $t+3k+r+2$ is integer. We have $|E(H)|=m+3n+r+s+1$. Therefore it is enough to prove the following
\begin{eqnarray*}
s    &   \leq  &   \frac{t+3k+r+2-\rho(m+3n+r+1)}{\rho}= \frac{t+3k+(1-\rho)r+1-\rho(m+3n)}{\rho} + \frac{1-\rho}{\rho}.
\end{eqnarray*}

\noindent We observe that the above inequality holds by the definition of $s$.

\noindent Assume now that $<H, k+1>$ is a yes-instance of PVC($\rho$). 
Then there exists a subset $C\subseteq V(H)$ of cardinality $k+1$ such that $C$ covers at least $\rho |E(H)|$ edges of $H$. We claim that there exists such a set $C$ with the additional property that $C$ contains the vertex $v$ (the central vertex of $K_{1,r}$ in the definition of the graph $H$) and also $|C\cap V(G)| = k$. To prove this claim first note that any vertex of degree one in $H$ covers exactly one edge; any vertex of $G$ covers at least three edges of $G'$; and also $k<n$. Hence by a modification of $C$ we may assume that $C$ contains no vertex of degree one and no vertex belonging to $P_s\cup \{w\}$. But from the other side, any vertex of $H$, other than $v$, covers no more than $n+2$ edges of $H$. So since $r\geq n+3$ if $v\notin C$ we may replace any vertex of $C$ by $v$. For this set $C$ we have $|C\cap V(G)|=k$.
Since $C$ covers at least $\rho |E(H)|$ edges of $H$, clearly $C\cap G$ covers at least $\lceil\rho |E(H)|\rceil -r-2$ edges of $G'$. Hence if we show that $\lceil\rho |E(H)|\rceil -r-2 \geq t+3k$ then we obtain that $<G', k, t+3k>$ is a yes-instance of PVCB. By the definition of $s$, there exists a value $\lambda$, $0\leq \lambda < 1$ such that
$$s=\frac{t+3k+(1-\rho)r+1-\rho(m+3n)}{\rho}-\lambda.$$\noindent
We have the following equivalent inequalities.
\begin{eqnarray*}
{} & {} &  \lceil\rho |E(H)|\rceil \geq t+3k+r+2  \\
{} &  \Leftrightarrow  &  \lceil\rho (m+3n+r+s+1)\rceil \geq t+3k+r+2  \\
{} &  \Leftrightarrow  &  \lceil\rho (m+3n+r+1)+t+3k+(1-\rho)r+1-\rho(m+3n)-\lambda\rho\rceil \geq t+3k+r+2  \\
{} &  \Leftrightarrow  &  \lceil\rho +t+3k+r+1-\lambda\rho\rceil \geq t+3k+r+2  \\
{} &  \Leftrightarrow  &  \lceil(1-\lambda)\rho\rceil \geq 1.
\end{eqnarray*}

\noindent Finally, the last inequality holds because $(1-\lambda)\rho>0$. This completes the proof.
\end{proof}

\noindent The NP-completeness result concerning PVCB is as follows.

\begin{thm}\label{PVCRHO}
PVC($\rho$) is NP-hard for bipartite graphs.
\end{thm}

\noindent \begin{proof}
For simplicity, denote the problem specified in the theorem by PVCB($\rho$). We make a reduction from PVCB to PVCB($\rho$). Assume that a bipartite graph $G=(V,E)$ on $n$ vertices with $m$ edges and two positive integers $t$ and $k$ are given. We may assume that $k\leq n$. From $<G,t,k>$ we obtain an instance of PVCB($\rho$) as follows.
Let $G'$ be obtained from $G$ as described in Lemma \ref{shakh}. Let also $H$ be the graph obtained from $G'$ as described in Lemma \ref{lem2}. Our desired instance for PVCB($\rho$) is $<H,k+1>$. In Lemma \ref{shakh} we showed that $<G, k, t>$ is a yes-instance of PVCB if and only if $<G', k, t+3k>$ is a yes-instance of PVCB. In Lemma \ref{lem2} we showed that $<G', k, t+3k>$ is a yes-instance of PVCB if and only if $<H, k+1>$ is a yes-instance of PVCB($\rho$). The proof completes.
\end{proof}

\noindent It was proved in \cite{WBG} that PVC is NP-hard for planar graphs and also for chordal graphs. Based on these results and the fact that transformations in Lemma \ref{shakh} and Lemma \ref{lem2} keeps the planarity and chordalness, hence by these lemmas and Theorem \ref{PVCRHO} we have the following remark.

\begin{remark}\label{rem1}
For any fixed $\rho$ with $0<\rho<1$, PVC($\rho$) is NP-hard even for the class of planar graphs and also the class of chordal graphs.
\end{remark}

\noindent We give the following family of bipartite graphs for which PVCB has a polynomial-time solution. Note that this family includes biregular bipartite graphs.

\begin{thm}
Let $G$ be a bipartite graph whose bipartite sets are $X$ and $Y$. Assume that the vertex degrees in $X$ are in decreasing form $d(x_1)\geq d(x_2)\geq\cdots\geq d(x_{|X|})$. Let $\delta(X)$ (resp. $\Delta(Y)$) be the minimum (resp. maximum) degree in $X$ (resp. $Y$). Assume that $\delta(X)\geq \Delta(Y)$. Then for any number $t$ $$P\beta_t(G)=\min\left\{k: {\sum}_{i=1}^k d(x_i)\geq t\right\}.$$\label{bireg}
\end{thm}

\noindent \begin{proof}
Set $k'=\min\left\{k: {\sum}_{i=1}^k d(x_i)\geq t\right\}$. It is obvious that $x_1,\ldots,x_{k'}$ covers at least $t$ edges of $G$. From the other hand, no set of size less than $k'$ covers $t$ edges of $G$. This completes the proof.
\end{proof}

\section{On the smallest monopolies with given average threshold}

\noindent In this section we first obtain a relationship between partial vertex covers in $G$ and $Smon_t(G)$. Next, we prove a hardness result concerning $Smon_t(G)$.

\begin{thm}\label{smon}
Let $G$ be any graph on $n$ vertices. Then for any positive number $t$ we have
$$Smon_t(G)=P\beta_{nt/2}(G).$$
\end{thm}

\noindent\begin{proof}
By the definition of $Smon_t(G)$, there corresponds a threshold assignment $\tau$ with ${\sum}_{v\in V(G)}\tau(v)\geq nt$ and a subset $M\subseteq V(G)$ such that $M$ is a $\tau$-monopoly of $G$ and $|M|=Smon_t(G)$. We define another threshold assignment $\varphi$ as follows. For each vertex $v\in V(G)\setminus M$ set $\varphi(v)= d_M(v)$ and for each vertex $v\in M$ set $\varphi(v)= d_G(v)$. It is clear that $M$ is a $\varphi$-monopoly of $G$. Also for each vertex $v\in V(G)$ we have $\tau(v)\leq \varphi(v)$. Hence ${\sum}_{v\in V(G)}\varphi(v)\geq nt$. In fact for any arbitrary subset $M$ we can define the threshold assignment $\varphi$ as above. We conclude that $Smon_t(G)$ equals to the smallest cardinality of a subset $M$ such that the threshold assignment $\varphi$ corresponding to $M$ satisfies ${\sum}_{v\in V(G)}\varphi(v)\geq nt$. Let $M_0$ be such that $Smon_t(G)=|M_0|$ and its corresponding $\varphi$ satisfies ${\sum}_{v\in V(G)}\varphi(v)\geq nt$. We note that the number of edges covered by $M_0$ is $\left({\sum}_{v\in V(G)}\varphi(v)\right)/2$. This shows that $P\beta_{nt/2}(G)\leq Smon_t(G)$. By a similar argument we obtain $Smon_t(G) \leq P\beta_{nt/2}(G)$. This completes the proof.
\end{proof}

\noindent The related complexity problem concerning $Smon$ is as follows, where $0<k<2$ is any fixed number.

\noindent {\bf Name: $Smon(k)$}\\
\noindent {\bf Instance:} A graph $G=(V,E)$ and a positive integer $d$.\\
\noindent {\bf Question:} Is there an assignment of thresholds $\tau$ to the vertices of $G$ with $n\bar{\tau}=\left\lceil nk\epsilon\right\rceil$ such that there exists a $\tau$-monopoly of size at most $d$?

\noindent Let $k$ be any fixed number with $0<k<2$. Theorem \ref{smon} shows that $Smon_{k\epsilon}(G)=P\beta_{km/2}(G)$. Hence the following is immediate corollary of Theorem \ref{PVCRHO} and Remark \ref{rem1}.

\begin{thm}
For any fixed number $k$ with $0<k<2$, $Smon(k)$ is NP-complete even restricted to the classes bipartite graphs, planar graphs and chordal graphs.
\end{thm}

\section{The smallest dynamic monopolies with given average\newline threshold}

\noindent In this section we first obtain a relationship between $Sdyn_t(G)$ and partial vertex covers in $G$. Next, we prove a hardness result concerning $Sdyn_t(G)$.

\begin{thm}\label{dynPVC}
Let $G$ be a graph of order $n$ and size $m$. Then $Sdyn_{t}(G)=P\beta_{nt-m}(G)$.
\end{thm}

\noindent \begin{proof}
As mentioned earlier in the paper, it was proved in \cite{Z2} that $Sdyn_t(G)=|G|-\max \{|H|: H\unlhd G, ~ |E(H)|\leq 2m-nt\}$.
Let $H_0$ maximize the related set in the latter equality. The proof of this result in \cite{Z2} shows that the complement of $H_0$, i.e. $H_0^c$ is a dynamic monopoly with the smallest size and with average threshold $t$. Obviously, $V(H_0^c)$ covers all the edges of $G$, except those edges with both endpoints in $H_0$. Hence, $H_0^c$ covers at least $m-(2m-nt)=nt-m$ edges. Then by the definition, $V(H_0^c)$ is an $(nt-m)$-partial vertex cover of $G$. It follows that $Sdyn_{t}(G)=P\beta_{nt-m}(G)$, as desired.
\end{proof}

\noindent The related complexity problem concerning $Sdynamo$ is as follows, where $1<k<2$ is any fixed number.

\noindent {\bf Name: $Sdynamo(k)$}\\
\noindent {\bf Instance:} A graph $G=(V,E)$ and a positive integer $d$.\\
\noindent {\bf Question:} Is there an assignment of thresholds $\tau$ to the vertices of $G$ with $n\bar{\tau}=\left\lceil nk\epsilon\right\rceil$ such that there exists a $\tau$-dynamic monopoly of size at most $d$?

\noindent Let $k$ be any fixed number with $1<k<2$. Theorem \ref{dynPVC} shows that $Sdyn_{k\epsilon}(G)=P\beta_{(k-1)m}(G)$. Hence the following is immediate corollary of Theorem \ref{PVCRHO} and Remark \ref{rem1}.

\begin{thm}
For any fixed number $k$ with $1<k<2$, $Sdynamo(k)$ is NP-complete even restricted to the three classes bipartite graphs, planar graphs and chordal graphs.
\end{thm}

\noindent It should be mentioned that when $k\leq 1$, then as proved in \cite{Z2}, $Sdyn_{k\epsilon}(G)=0$. It was proved in \cite{CS} that $PVC$ admits a polynomial-time algorithm for trees. Therefore $Smon(k)$ and $Sdynamo(k)$ have polynomial-time solutions for trees. These problems can also be solved in polynomial-time for biregular bipartite graphs by Theorem \ref{bireg}.

\end{document}